\documentclass[12pt,a4paper]{amsart}
\usepackage{amssymb}
\usepackage{mathrsfs}

\headsep=1cm \numberwithin{equation}{section}
\newtheorem{thm}{Theorem}[section]
\newtheorem{cor}[thm]{Corollary}
\newtheorem{lem}[thm]{Lemma}

\theoremstyle{definition}

\theoremstyle{remark}

\numberwithin{equation}{section}

\newcommand\Ass{\operatorname{Ass}}

\newcommand\Hom{\operatorname{Hom}}
\newcommand\Ext{\operatorname{Ext}}

\newcommand\m{\operatorname{\frak m}}

\begin{document}
\title{A new characterization of Cohen-Macaulay rings}
\author[Kamal Bahmanpour and Reza Naghipour]{Kamal Bahmanpour and Reza Naghipour$^*$}
\address{Department of Mathematics, Ardabil branch, Islamic Azad University, P.O. Box 5614633167,
Ardabil, Iran; and School of Mathematic, Institute for studies in Theoretical
Physics and Mathematics (IPM), P.O. Box 19395-5746, Tehran, Iran.} \email{bahmanpour.k@gmail.com}

\address{Department of Mathematics, University of Tabriz, Tabriz, Iran;
and School of Mathematic, Institute for studies in Theoretical
Physics and Mathematics (IPM), P.O. Box 19395-5746, Tehran, Iran.}
\email{naghipour@ipm.ir} \email {naghipour@tabrizu.ac.ir}

\thanks{
This research  was been in part supported by a grant from IPM (No. 91130022)\\
$^*$Corresponding author: e-mail: {\it naghipour@ipm.ir} (Reza
Naghipour)}%

\begin{abstract}
The purpose of this article is to provide a new characterization of Cohen-Macaulay local rings. As a consequence
we deduce that a local (Noetherian) ring $R$ is Gorenstein if and only if every parameter ideal of $R$ is irreducible.
\end{abstract}
\maketitle
\section{Introduction}
Let $(R,\m)$ denote a commutative Noetherian local ring.  It is well-known that $R$ is Gorenstein if and only if $R$ is
Cohen-Macaulay and some ideal generated by a system of parameters (s.o.p) of $R$ is irreducible. Perhaps less
widely known is a nice result of D.G. Northcott and D. Rees which states that if every ideal generated by a s.o.p
of $R$ (henceforth a {\it parameter ideal}) is irreducible, then $R$ is Cohen-Macaulay (see  \cite[Theorem 1]{NR}). Thus,
$R$ is Gorenstein  if and only if every parameter ideal is irreducible.
The purpose of this note is to establish a new characterization of Cohen-Macaulayness of $R$. Our main theorem is the following:

\begin{thm}
Let $(R, \frak m)$ be a  local (Noetherian) ring. Then the following conditions are equivalent:

{\rm(i)} $R$ is a Cohen-Macaulay ring.

{\rm(ii)} For any parameter ideals $\frak a$ and $\frak b$ of $R$ with $\frak b\subseteq \frak a$, we have $$\Hom_R(R/\frak a,R/\frak b)\cong R/\frak a.$$
\end{thm}

As a consequence we obtain that a local (Noetherian) ring $R$ is Gorenstein if and only if every parameter ideal of $R$ is irreducible.

One of our tools for proving Theorem 1.1 is the following:

\begin{lem}
Let $(R, \frak m)$ be a  local (Noetherian) ring.

{\rm(i)} If $x\in \m$, then there exists a non-negative integer $n$ such that $(x^n)\cap \Gamma_{\m}(R)=0$.

{\rm(ii)} If $\dim\,R=d\geq 1 $ and $\Gamma_{\m}(R)\neq 0$, then there exists a system of parameters $x_1,\dots,x_d$ of $R$ such that $\m\in \Ass_R\,R/(x_1,\dots,x_i)$ for all $1\leq i \leq d$.
\end{lem}

Throughout this article, $R$ will always be a commutative Noetherian local ring with maximal ideal $\m$ and $\dim R=d$.
 For any $R$-module $L$, the submodule $\bigcup_{n\geq 1}(0:_L\m^n)$ of $L$ is denoted by $\Gamma_{\m}(L)$.
 Also, the socle of $L$ is defined to be  $(0:_L\m)=\{x\in L|\, \m x=0\}$, and is denoted by ${\rm Soc}_R\, L$.
  For any unexplained notation and terminology we refer the reader to \cite{BH} and \cite{Mat}.

\section{Proof of the main theorem}

In order our main result, we prepare the following lemma.

\begin{lem}
{\rm(i)} If $x\in \m$,  then there exists an integer $n\geq1$ such that $(x^n)\cap \Gamma_{\m}(R)=0$.

{\rm(ii)} If $\dim R=d\geq 1$ and $\m \in \Ass_RR$, then there exists a system of parameters $a_1,\dots, a_d$ of $R$ such that $\m\in \Ass_RR/(a_1,\dots, a_i)$ for all $1\leq i \leq d$.
\end{lem}
\proof (i) Since $\Gamma_{\m}(R)$ is an Artinian $R$-module, the descending sequence of the submodules of $\Gamma_{\m}(R)$,$$(x)\cap \Gamma_{\m}(R)\supseteq (x^2)\cap \Gamma_{\m}(R)\supseteq \cdots$$
is stationary, i.e., there exists a positive integer $n$ such that $$(x^{n+i})\cap \Gamma_{\m}(R)=(x^n)\cap \Gamma_{\m}(R),$$ for all integers $i \geq 1$. Now, we have $$(x^n)\cap \Gamma_{\m}(R)\subseteq \bigcap_{i\geq n}(x^{i})\cap \Gamma_{\m}(R)\subseteq \bigcap_{i\geq 1}(x^{i}),$$and so by Krull's Intersection Theorem $(x^n)\cap \Gamma_{\m}(R)=0$, as required.

(ii) Let $I=\Gamma_{\m}(R)$. Since $\m \in \Ass_R\,R$, it follows that $I\neq0$. Now let $b_1,\dots, b_d$ be a system of parameters of $R$. Then in view of (i) there exists an integer $n_1\geq 1$ such that $(b_1^{n_1})\cap I=0$. Set $a_1:=b_1^{n_1}$. Then $I+(a_1)/(a_1)\neq 0$ and $I+(a_1)/(a_1)\subseteq \Gamma_{\m}(R/(a_1))$. It is clear that $a_1$ is a part of a system of parameters of $R$. Let $R_1:=R/(a_1)$ and $\m_1:=\m/(a_1)$. Then $\Gamma_{\m}(R_1)=\Gamma_{\m_1}(R_1)\neq 0$ and so $\m_1\in \Ass_{R_1}R_1$. Moreover, the elements $b_2+(a_1),\dots, b_d+(a_1)$ is a system of parameters of $R_1$. (Note that $\dim \,R/(a_1)=d-1$.) Therefore, there is an integer $n_2\geq 1$ such that $b_2^{n_2}R_1\cap\Gamma_{\m_1}(R_1)= 0$. Now, we put $a_2:=b_2^{n_2}$. Then $\Gamma_{\m}(R/(a_1,a_2))\neq 0$ and $a_1,a_2$ is a part of a system of parameters of $R$. Continue in this way: there exist positive integers $n_1,\dots, n_d$ such that $a_1=b_1^{n_1},\dots, a_d=b_d^{n_d}$ is a system of parameters of $R$ and $\m\in \Ass_R\,R/(a_1,\dots, a_i)$ for all $1\leq i \leq d$.
\qed\\

\begin{cor} Suppose that for every system of parameters $a_1,\dots, a_d$ of $R$, the element $a_d$ is an $R/(a_1,\dots, a_{d-1})$-regular. Then $\m\not\in \Ass_RR$.
\end{cor}
\proof Suppose that $\m\in \Ass_RR$ and look for a contradiction. To do this, in view of Lemma 2.1, there exists a system of parameters $x_1,\dots, x_d$ of $R$ such that $\m\in \Ass_R\,R/(x_1,\dots, x_{d-1})$, and so $x_d$ is not an $R/(x_1,\dots, x_{d-1})$-regular element. With this contradiction the proof is complete.\qed\\

Now, we are ready to state and prove the main theorem.

\begin{thm}
The following conditions are equivalent:

{\rm(i)} $R$ is a Cohen-Macaulay ring.

{\rm(ii)} For any parameter ideals $\frak a$ and $\frak b$ of $R$ with $\frak b\subseteq \frak a$, we have
$$\Hom_R(R/\frak a,R/\frak b)\cong R/\frak a.$$
\end{thm}
\proof First we show (i)$\Rightarrow$(ii). To this end, let $\frak a=(a_1,\dots, a_d)$ and $\frak b=(b_1,\dots, b_d)$ be two parameter ideals of $R$ with $\frak b \subseteq \frak a$.
Then, as $R$ is a Cohen-Macaulay ring of dimension $d$, it follows that $\{a_i\}_{i=1}^d$ and  $\{b_i\}_{i=1}^d$ are two maximal $R$-sequences in $\frak a$. Consequently, in view of Rees' theorem (see \cite[Lemma 1.2.4]{BH}), we have $$\Ext^d_R(R/\frak a,R)\cong \Hom_{R/\frak a}(R/\frak a,R/\frak a) \cong R/\frak a,\,\,\,\,\,{\rm and}$$
$$\Ext^d_R(R/\frak a,R)\cong \Hom_{R/\frak b}(R/\frak a,R/\frak b) = \Hom_{R}(R/\frak a,R/\frak b).$$
Therefore, $\Hom_R(R/\frak a,R/\frak b)\cong R/\frak a,$ as required.

In order to prove (ii)$\Rightarrow$(i), we use induction on $d$. When $d=1$, there exists $a\in R$ such that $\sqrt{(a)}=\m$. Then, in view of Lemma 2.1 there is a positive integer $n$ such that $(a^n)\cap I=0$, where $I=\Gamma_{\m}(R)$. Moreover, there exists an integer $k$ such that $\m^kI=0$, and so $a^kI=0$. Let $r_1=\max\{n,k\}$ and $b=a^{r_1}$. Then $\sqrt{(b)}=\m$, $\Gamma_{(b)}(R)=I$, $bI=0$ and $(b)\cap I=0$. On the other hand, in view of \cite[Proposition 4.7.13]{No}, there exists an integer $r_2\geq 1$ such that for all integers $t\geq r_2$,
$$(b^t):_Rb=b^{t-r_2}((b^{r_2}):_Rb)+(0:_Rb).$$
As $(0:_Rb)=I$, it follows that $$(b^t):_Rb=b^{t-r_2}((b^{r_2}):_Rb)+I.$$
Now, let $s\geq r_2+1$ be an integer. Then in view of assumption (ii) we have $$R/(b)\cong \Hom_R(R/(b),R/(b^s))\cong \frac{(b^s):_Rb}{(b^s)}=\frac{b^{s-r_2}((b^{r_2}):_Rb)+I}{(b^s)}.$$
Next, we show that $$b^{s-r_2}((b^{r_2}):_Rb)\cap(I+(b^s))=(b^s).$$
To do this, suppose $c\in I$ and $\delta b^{s-r_2}\in b^{s-r_2}((b^{r_2}):_Rb)$ such that $\delta b^{s-r_2}=c+ub^s$, where $\delta \in (b^{r_2}):_Rb$ and $u\in R$. Then $$c=\delta b^{s-r_2}- ub^s\in (b)\cap I,$$
and so $c=0$. Thus $\delta b^{s-r_2}=ub^s\in (b^s)$. Consequently, $$\frac{R}{(b)}\cong \frac{b^{s-r_2}((b^{r_2}):_Rb)}{(b^s)}\oplus \frac{I+(b^s)}{(b^s)}.$$
Since $R/(b)$ is a local ring, it follows that $$b^{s-r_2}((b^{r_2}):_Rb)=(b^s)\,\,\,\,\,{\rm or}\,\,\,\,\,I+(b^s)=(b^s).$$ If $b^{s-r_2}((b^{r_2}):_Rb)=(b^s)$, then as $b^{s-1}\in b^{s-r_2}((b^{r_2}):_Rb)$, it follows that $(b^{s-1})=(b^s)$. Hence in view of NAK's Lemma, $b^{s-1}=0$, and so $\m=\sqrt{(0)}$, which is a contradiction. Therefore, $I\subseteq (b^s)$, and so $I=I\cap(b^s)\subseteq I\cap (b)=0$. Thus $I=0$, and so $\m\not\in \Ass_RR$. Since $\dim R=1$, it follows that $R$ is a Cohen-Macaulay ring.

Assume, inductively, that $d\geq 2$ and that the result has been proved for $d-1$. Let $x_1,\dots, x_d$ be an arbitrary system of parameters for $R$. For each integer $n\geq 1$, set $$\frak a:=(x_1,\dots, x_d) \,\,\,\,\,{\rm and}\,\,\,\,\, \frak b:=(x_1,\dots, x_{d-1},x_d^n).$$ Then $$\Hom_R(R/\frak a,R/\frak b)\cong R/\frak a.$$ Let $S=R/(x_1,\dots, x_{d-1})$. Then $$\Hom_S(S/x_dS,S/x_d^nS)\cong \Hom_R(R/\frak a,R/\frak b)\cong R/\frak a \cong S/x_dS.$$ Now, if $\overline{y}:=y+(x_1,\dots, x_{d-1})$ is a system of parameters of $S$, then $x_1,\dots, x_{d-1},y$ is a system of parameters of $R$, (note that $\dim S=1$). Hence$$\Hom_S(S/(\overline{y}), S/(\overline{y})^n)\cong S/(\overline{y}),$$
and so applying the method used above it follows that $S$ is a Cohen-Macaulay ring. Thus $x_d$ is an $R/(x_1,\dots, x_{d-1})$-regular element, and hence by Corollary 2.2, $\m\not\in \Ass_R\,R$. Accordingly, there exists $\eta_1\in \m$ such that $\eta_1\not\in Z_R(R)$. Therefore, $\eta_1$ is a part of a system of parameters of $R$. Now, let $T=R/(\eta_1)$ and $$\frak c=(\bar \eta_2,\dots, \bar\eta_d),\,\,\,\,\,{\rm and}\,\,\,\,\,\frak d=(\bar \xi_2,\dots,\bar\xi_d)$$ be two parameter ideals of $T$ such that $\frak d \subseteq \frak c$, where $\bar \eta_i= \eta_i+ (\eta_1)$ and $\bar \xi_i= \xi_i+ (\eta_1)$ for every $i=2, \dots, d$. Then $$\frak a'=(\eta_1,\dots,\eta_d)\,\,\,\,\,{\rm and}\,\,\,\,\,\frak b'=(\eta_1,\xi_2,\dots,\xi_d)$$ are two parameter ideals of $R$. Moreover, we have $$\frak b'\subseteq \frak a'\,\,\,\,\,{\rm and}\,\,\,\,\,\Hom_R(R/\frak a',R/\frak b')\cong R/\frak a'.$$
Therefore $$\Hom_T(T/\frak c,T/\frak d)\cong \Hom_R(R/\frak a',R/\frak b')\cong R/\frak a'\cong T/\frak c.$$ As, $\dim\,T=d-1$, it follows from the inductive hypothesis that $T$ is a Cohen-Macaulay ring. Hence $R$ is a Cohen-Macaulay ring, as required. \qed\\

\begin{cor}
The following conditions are equivalent:

{\rm(i)} $R$ is a Cohen-Macaulay ring.

{\rm(ii)} For any parameter ideals $\frak a$ and $\frak b$ of $R$ with $\frak b\subseteq \frak a$ the $R$-module
$\Hom_R(R/\frak a,R/\frak b)$ is cyclic.
\end{cor}
\proof The result follows from the proof of Theorem 2.3.\qed\\

\begin{cor}
The following conditions are equivalent:

{\rm(i)} $R$ is a Gorenstein ring.

{\rm(ii)} For any parameter ideals $\frak a$ and $\frak b$ of $R$ with $\frak b\subseteq \frak a$,
$\Hom_R(R/\frak a,R/\frak b)\cong R/\frak a$ and there exists an irreducible parameter ideal.

{\rm(iii)} For any parameter ideals $\frak a$ and $\frak b$ of $R$ with $\frak b\subseteq \frak a$ the $R$-module
$\Hom_R(R/\frak a,R/\frak b)$ is cyclic and there exists an irreducible parameter ideal.
\end{cor}
\proof The result follows from the proof of Theorem 2.3 and \cite[Theorem 18.1]{Mat}.\qed\\

The final result shows  that a local (Noetherian) ring $R$ is Gorenstein if and only if every parameter ideal of $R$ is irreducible.
\begin{cor}
$R$ is Gorenstein if and only if  every parameter ideal of $R$ is irreducible.
\end{cor}
\proof In view of Corollary 2.5, it is enough to show that for any parameter ideals $\frak a$ and $\frak b$ of $R$ with $\frak b\subseteq \frak a$ we have
$\Hom(R/ \frak a, R/\frak b)\cong R/\frak a$. To do this end, since  $\frak a$ and $\frak b$ are irreducible it follows that
$$\dim_{R/\m}\, {\rm Soc}_R\, R/\frak a=1=\dim_{R/\m}\,{\rm Soc}_R\, R/\frak b.$$  In particular, $R/\frak a$ and $R/\frak b$ are two Gorenstein local rings of dimension zero. Hence $E_{R/\frak a}(R/\m)\cong R/\frak a$ and $E_{R/\frak b}(R/\m)\cong R/\frak b$. Therefore,

\begin{eqnarray*}
\Hom_R(R/\frak a,R/\frak b)&\cong&\Hom_{R/\frak b}(R/\frak a,R/\frak b)\\ &\cong& \Hom_{R/\frak b}(R/\frak a,E_{R/\frak b}(R/\m))\\ &\cong& E_{R/\frak a}(R/\m)\\ &\cong& R/\frak a.\end{eqnarray*}
This completes the proof.

\begin{center}
{\bf Acknowledgments}
\end{center}
We would like to thank from School of Mathematics, Institute for Research in Fundamental
Sciences (IPM), for its financial support.


\end{document}